\documentclass[12pt,a4paper]{article}

\usepackage[utf8]{inputenc}
\usepackage[T1]{fontenc}
\usepackage{lmodern}

\usepackage{amsmath, amssymb, amsthm}
\usepackage{authblk}
\usepackage{float}
\usepackage{array}
\usepackage{booktabs}
\usepackage{graphicx}      
\usepackage{subcaption}    
\usepackage{caption}       
\usepackage{tikz}
\usepackage[a4paper, margin=1in]{geometry}

\usepackage[hidelinks]{hyperref}
\usepackage{longtable}

\usepackage{sectsty}
\allsectionsfont{\normalfont\scshape}

\usepackage{algorithm}
\usepackage{algpseudocode}

\usepackage{float} 
\usepackage[T1]{fontenc}

\theoremstyle{remark} 
\newtheorem{remark}{Remark}
\theoremstyle{example} 
\newtheorem{example}{Example}

\title{Extending Hridaya Kolam to Even-Ordered Dot Patterns and Their Applications}

\author[1]{Suvra Kanti Chakraborty}
\author[2]{Atanu Manna\thanks{Corresponding author: atanuiitkgp86@gmail.com; atanu.manna@iict.ac.in}}
\affil[1]{\small Ramakrishna Mission Vidyamandira, Belur Math- 711202, West Bengal, India}
\affil[2]{\small Indian Institute of Carpet Technology, Chauri Road, Bhadohi--221401, Uttar Pradesh, India}
\date{} 
\begin{document}
\maketitle

\begin{abstract}
\noindent This study extends the mathematical framework of Hridaya Kolam patterns by applying modular arithmetic to even-ordered dot arrangements with arm counts co-prime to the number of dots. We analyze the resulting cyclic sequences that correspond to Eulerian circuits, enabling continuous single-stroke kolam designs beyond the classical odd-ordered cases. Our method provides explicit algorithms for constructing these intricate patterns, unveiling new symmetries and structural properties. Elevating this traditional floor art, we translate these mathematically grounded motifs into striking designs, showcasing their beauty and complexity in contemporary dari art in the carpet and textile sectors.
\end{abstract}

$\\$\textbf{Keywords:} Kolam; Modular arithmetic; Sequence; Directed graph.

\section{Introduction}\label{secint}

Kolam, also known in various regions as Kamalam, is a revered traditional floor art form deeply embedded in the cultural and spiritual fabric of India. Predominantly practiced in southern states such as Tamil Nadu, Kerala, Karnataka, and parts of Andhra Pradesh and Telangana, Kolam patterns are meticulously crafted using rice flour or limestone powder. These intricate, symmetrical, and often geometric designs are drawn daily or during auspicious occasions to invoke prosperity, peace, and positive energy within homes and communal spaces. Kolam is part of a diverse folk art tradition that varies widely across regions and is known by distinct names such as Muggulu or Muggu in Andhra Pradesh and Telangana, Alpana in West Bengal, Puvidal in Kerala, Chowk-poorana in Punjab and northern states, Mandana in Rajasthan, Rangoli in Maharashtra and Karnataka, and Sanjhi in Uttar Pradesh. 
Despite stylistic differences, these designs function uniformly as profound visual expressions of devotion, cultural identity, and aesthetic harmony. The use of natural materials and the underlying symbolic geometry have drawn the interest of artists and researchers in anthropology, ethnography, design theory, and computational geometry. When colored powders replace traditional white media, Kolam transforms into the more widely recognized Rangoli or Rangavalli.\\


\noindent The pioneering work by Siromoney \cite{SIROMONEY78}, and Siromoney and Chandrasekaran \cite{SIROCHANDRA} presented at the 1986 ISI Calcutta conference systematically analyzed the Hridaya Kamalam Kolam—a distinct South Indian design constructed via dots arranged radially on multiple arms. Their investigation primarily focused on configurations with an odd number of dots per arm, where the arm count is coprime to the dot count. 

\begin{figure}[H]
    \centering
    \begin{subfigure}{0.30\textwidth}
        \centering
        \includegraphics[width=1\linewidth]{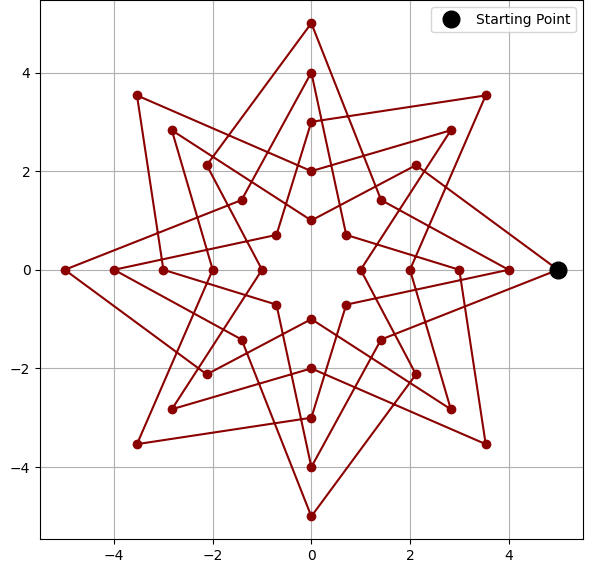}
        \caption{Dot-Arm structure}
        \label{fig:dot_arm_structure}
    \end{subfigure}\hfill
    \begin{subfigure}{0.30\textwidth}
        \centering
        \includegraphics[width=\linewidth]{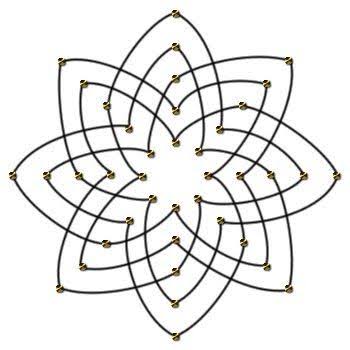}
        \caption{Curvilinear design}
        \label{fig:kolam_digital}
    \end{subfigure}\hfill
    \begin{subfigure}{0.30\textwidth}
        \centering
        \includegraphics[width=\linewidth]{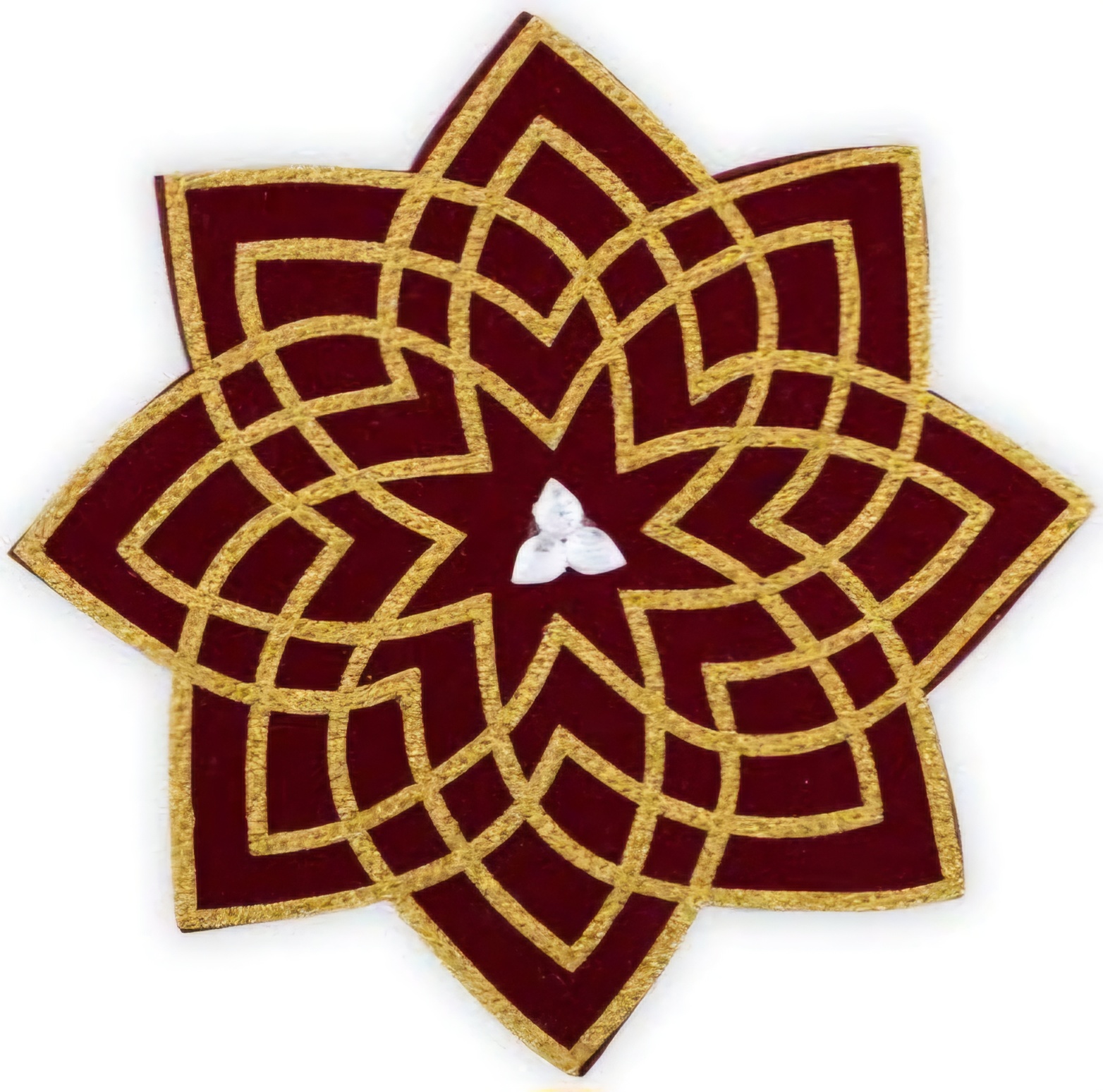}
        \caption{Kolam artifact}
        \label{fig:kolam_real}
    \end{subfigure}
    \caption{Hridaya Kolam ($m=5~ dots~(odd), n=8~ arms):  5\rightarrow3\rightarrow1\rightarrow4 \rightarrow 2$}
    \label{fig:kolam_full_view}
\end{figure}

\noindent Significant contributions relevant to this domain include the formal language approach to Kolam generation by Ascher \cite{ASCHER}, the Fibonacci sequence-based pattern constructions by Naranan and collaborators \cite{NARANAN, NARANAN2, NARANAN3, NARANAN4}, and the numerical and diagrammatic analyses by Yanagisawa and Nagata \cite{YANAGISAWA}.
Robinson’s ‘Pasting Scheme’ \cite{ROBINSON} and Nagata’s digitalization studies \cite{NAGATA, NAGATA2} further advanced computational treatments of these patterns. Mathematical insights connecting Kolam design with graph theory and algebra have been explored by Thirumuthi and Simic-Muller \cite{THIRUMURTHY}, and Sarin \cite{SARIN}.\\\\
\noindent Recently in 2023, Srinivasan’s algorithmic approach to Hridaya and Aishwarya Kolams \cite{SRINIVASAN} has facilitated the efficient generation and analysis of these patterns, offering deeper insights into their underlying structural principles. Specifically, this work of Srinivasan focuses on configurations of Hridaya Kolams, where the number of dots arranged in a circular layout, denoted by \( m \), is \emph{odd}, and the arm count \( n \) satisfies \( \gcd(m, n) = 1 \). Moreover, Srinivasan observed that a mathematical sequence underlies Kolam patterns and provided an algorithm for generating such sequences. Additionally, he introduced a generalized sequence to produce various novel Hridaya Kolam designs. However, his analysis is restricted to cases where the number of dots \( m \) is odd; configurations with even values of \( m \) remain unexplored. This naturally leads to the following questions:\\[.7em]
\noindent \textsf{Q(a) Can Hridaya Kolam patterns be generated for configurations with an even number of dots (\( m \))?}\\[.2em] and\\[.2em]
\noindent \textsf{Q(b) Is there a novel approach to generate the underlying sequences for these patterns?} \\\\
\noindent Therefore, the main objective of this paper is to address the questions posed in \textsf{Q(a)} and \textsf{Q(b)} above. To this end, we employ the framework of congruence relations and modular arithmetic to develop methods for generating the underlying sequences, thereby providing an answer to \textsf{Q(b)}. Furthermore, this approach successfully produces Hridaya Kolam patterns with an even number of dots, addressing \textsf{Q(a)} as well. Remarkably, this approach also naturally covers the odd-numbered configurations, thereby offering a unified method applicable to both even and odd cases.\\\\
\noindent The paper is organized as follows. Section~2 presents the connection between modular arithmetic and the Hridaya Kolam. Section~3 details the main results of the paper, including the development of a systematic algorithm for generating generalized Kolam patterns. We analyze the resulting visual motifs, exploring the distinct geometric symmetries and aesthetic variations introduced by these new configurations. Section~4 demonstrates the practical application of these mathematically derived designs in the carpet industry, showcasing their translation into digitally crafted textile products such as Daries. Finally, Section~5 summarizes our findings, reflects on their implications for both cultural heritage and modern design, and suggests directions for future research.

\section{Modular Arithmetic and Hridaya Kolam}

\noindent 
In this section, we discuss how Hridaya Kolam is greatly connected with modular arithmetic. To begin with, we recall first the notion of modular arithmetic, in particular congruence relation as below. Given \( m \) as the number of dots per arm and \( n \) as the number of arms, define the sequence:
\[
a_k = 
\begin{cases}
m & \text{if } k = 0, \\
(k \cdot n) \bmod m & \text{if } k = 1, 2, \dots, m-1.
\end{cases}
\]
This congruence relation produces a single-stroke Kolam path through a permutation of \( \{1,2,\ldots,m\} \) if \( \gcd(m,n) = 1 \). For example, if one chooses \( m = 5 \) and \( n=8 \), then one such permutation is:
\(
\{0, 3, 1, 4, 2\} \equiv \{5, 3, 1, 4, 2\}
\).
\begin{figure}[H]
\centering
\begin{tikzpicture}[scale=2.5, every node/.style={font=\small}]
    \def \n {8} 
    \def \m {5} 
    \def \gap {0.25} 

    \foreach \i in {0,...,7} {
        \foreach \j in {1,...,5} {
            \filldraw ({\j*\gap*cos(360/\n*\i)}, {\j*\gap*sin(360/\n*\i)}) circle (0.5pt);
        }
        \draw[gray, thick, ->] (0,0) -- ({(\m+0.5)*\gap*cos(360/\n*\i)}, {(\m+0.5)*\gap*sin(360/\n*\i)});
    }

    \filldraw (0,0) circle (0.6pt);
\end{tikzpicture}
\hspace{2em}
\includegraphics[width=0.45\textwidth]{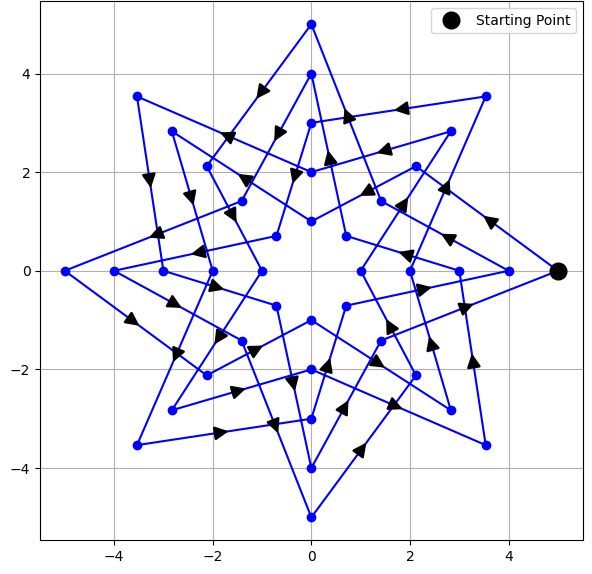}
\caption{Left: Radial dot configuration with \( n = 8 \) arms and \( m = 5 \) dots per arm. \quad Right: Corresponding Hridaya Kolam for the cyclic sequence $5\rightarrow3\rightarrow1\rightarrow4\rightarrow 2$ }
\end{figure}
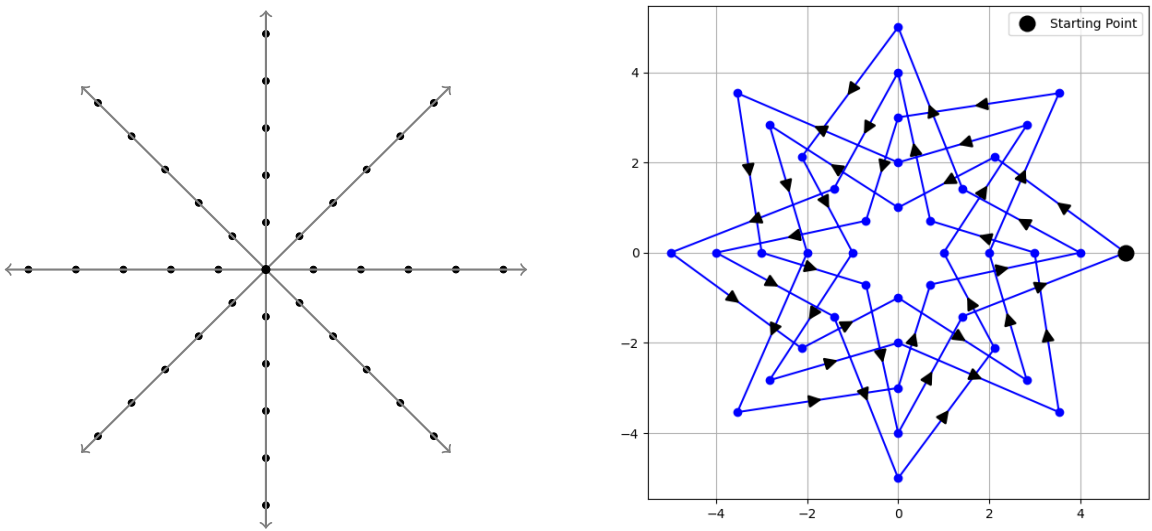

\noindent From an abstract perspective, this construction is grounded in the structure of the cyclic group \( \mathbb{Z}_m \). When \( \gcd(m, n) = 1 \), the element \( n \) serves as a generator, producing a complete residue system modulo \( m \). This guarantees that the modular sequence \( \{kn \bmod m\} \) cycles through all distinct residues, thus encoding the underlying Kolam pattern within the algebraic framework of modular arithmetic.

\begin{remark}
    The traditional Hridaya Kolam corresponds to the cyclic sequence 
    $1 \rightarrow 3 \rightarrow 5 \rightarrow 2 \rightarrow 4$, 
    which represents the same closed loop as described above, traced in a clockwise direction starting from $1$. 
    Since the loop is closed, the direction of traversal, whether clockwise or anticlockwise, does not affect the resulting Kolam pattern.
\end{remark}

\section{Hridaya Kolam for even number of dots}

Let us suppose that \( m \) is even and $n$ is such that \( \gcd(m, n)=1 \).  The Hridaya Kolam dot arrangement can be systematically constructed using the sequences generated by modular arithmetic as described in the previous section. These sequences determine radial positions for dots distributed evenly over \( n \) arms.\\

\noindent We begin with the above \( m \)-term sequence and repeat it \( n \) times to populate an \( m \times n \) matrix row-wise, where each column corresponds to a radial arm and each row represents one concentric layer of dots placed around the center. To visualize this, think of the matrix as a one-dimensional array formed by appending its rows consecutively. We insert the sequence repeatedly into this array in the following two cases:

\begin{enumerate}
    \item if \( m > n \), a single sequence does not fill a row completely, so the next sequence continues in the same row.
    \item if \( m < n \), the sequence overflows to the next row once a row is filled.
\end{enumerate}
\noindent   
Repetition of the sequence \( n \) times in this fashion allows us to construct the entire \( m \times n \) matrix with the following properties:

\begin{enumerate}
    \item Each row is a cyclic right-shift of the previous row.
    \item Each column (or “arm”) contains all values from the set \( \{1, 2, \dots, m\} \) exactly once.
    \item All \( m \cdot n \) positions in the matrix are uniquely filled without overlap or repetition.
\end{enumerate}
 \noindent For the sake of completeness, we now provide some examples of such matrices constructed with an even number of dots.
 
 \begin{example}   
 Choose \( m = 4 \), \( n = 3 \). In this case, $m>n$. By modular arithmetic concept, we have $a_k=3k \bmod 4$ for $k=0, 1, 2, 3$ with $a_0=4$. Hence the resulting sequence is $4\rightarrow3\rightarrow2\rightarrow1$, and the corresponding \( 4 \times 3 \) matrix is presented as below:
\[
\begin{array}{|c|c|c|}
\hline
\textbf{Arm 1} & \textbf{Arm 2} & \textbf{Arm 3} \\
\hline
4 & 3 & 2 \\
\hline
1 & 4 & 3 \\
\hline
2 & 1 & 4 \\
\hline
3 & 2 & 1 \\
\hline
\end{array}
\]
By closing the loop, the corresponding Hridaya Kolam is shown at Figure 3.
\begin{figure}[H]
\centering
\includegraphics[width=0.45\textwidth]{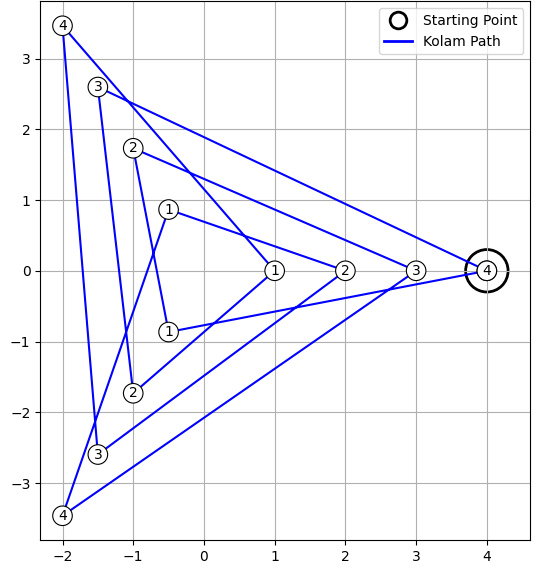}
\caption{Hridaya Kolam configuration for even \( m=4, n=3 \).}
\end{figure}
\end{example}

\begin{example}
Let us now consider an example for the case $m<n$. Choose \( m = 4 \), \( n = 5 \). Then, applying the same modular arithmetic technique, the corresponding sequence will be $4\rightarrow1\rightarrow2\rightarrow3$. Hence the resulting matrix of order \( 4 \times 5 \) takes the following form:
\[
\begin{array}{|c|c|c|c|c|}
\hline
\textbf{Arm 1} & \textbf{Arm 2} & \textbf{Arm 3} & \textbf{Arm 4} & \textbf{Arm 5} \\
\hline
4 & 1 & 2 & 3 & 4 \\
\hline
1 & 2 & 3 & 4 & 1 \\
\hline
2 & 3 & 4 & 1 & 2 \\
\hline
3 & 4 & 1 & 2 & 3 \\
\hline
\end{array}
\] 
and the corresponding Hridaya Kolam is depicted in the following Figure 4.
\begin{figure}[H]
\centering
\includegraphics[width=0.45\textwidth]{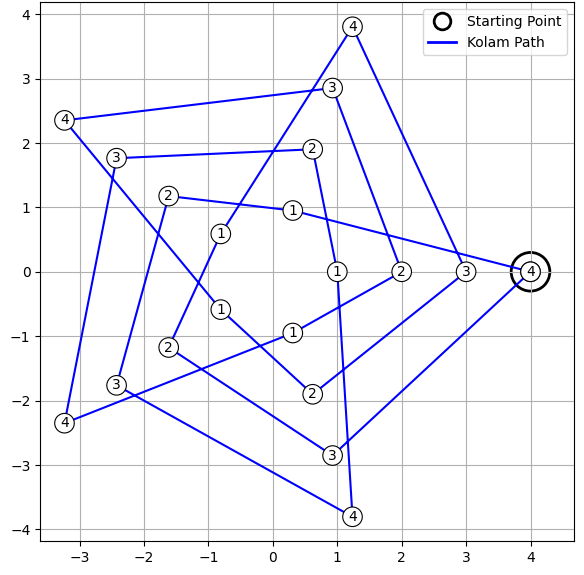}
\caption{Hridaya Kolam configuration for even \( m=4, n=5 \).}
\end{figure}
\end{example}
\noindent
The above two examples, in fact the above two Hridaya Kolams, gives us the following observations. First, each row represents one radial layer of dots placed around the center, starting from angle \( 0 \) and increasing by \( 2\pi/n \) for each arm. Secondly, the dots are placed at radial distances defined by the matrix entries. Since each arm (column) receives all values \( \{1, 2, \dots, m\} \) without repetition, so there is no dot duplicated on an arm. As there are \( m \) rows and \( n \) columns, all \( m\cdot n \) positions are used exactly once, guaranteeing non overlapping dot placements.\\

\noindent This structure ensures a visually elegant and mathematically robust Hridaya Kolam design, where each arm and layer are fully utilized. The traversal of adjacent dots results in a closed circuit forming an Eulerian path through the kolam graph.

\bigskip

\subsection{Modular Arithmetic-Based Generator Sequences}
In this section, we construct generating sequences for even values of \( m \) and integers \( n \) such that \( \gcd(m, n) = 1 \). The final term, identical to the initial term, is included deliberately to reinforce the looping behavior inherent to these sequences, ensuring that the design returns to its starting point and forms a complete traversal over the dot pattern.

\begin{longtable}{|c|c|p{11cm}|}
\caption{Sequences generated for various $m$ and $n$} \label{tab:modular_sequences} \\
\hline
\( m \) & \( n \) & Sequence (cycle) \\
\hline
\endfirsthead

\multicolumn{3}{c}%
{{\bfseries \tablename\ \thetable{} -- continued from previous page}} \\
\hline
\( m \) & \( n \) & Sequence (cycle) \\
\hline
\endhead

\hline \multicolumn{3}{r}{{Continued on next page}} \\
\endfoot

\hline
\endlastfoot

2 & 3, 5, 7, 9, 11, 13 & \(2\rightarrow1\rightarrow2\) \\
\hline
4 & 3, 7, 13 & \(4\rightarrow3\rightarrow 2\rightarrow 1 \rightarrow4\) \\
  & 5, 9, 11 & \(4\rightarrow 1\rightarrow 2\rightarrow 3\rightarrow 4\) \\
\hline
6 & 5, 11 & \(6\rightarrow5\rightarrow4\rightarrow3\rightarrow2\rightarrow 1\rightarrow6\) \\
  & 7, 13 & \(6\rightarrow1\rightarrow2\rightarrow3\rightarrow4\rightarrow5\rightarrow6\) \\
\hline
8 & 3, 11 & \(8\rightarrow3\rightarrow6\rightarrow1\rightarrow4\rightarrow7\rightarrow2\rightarrow5\rightarrow8\) \\
  & 5, 13 & \(8\rightarrow5\rightarrow2\rightarrow7\rightarrow4\rightarrow 1\rightarrow6\rightarrow3\rightarrow8\) \\
  & 7 & \(8\rightarrow7\rightarrow6\rightarrow5\rightarrow4\rightarrow 3\rightarrow2\rightarrow1\rightarrow8\) \\
  & 9 & \(8\rightarrow1\rightarrow2\rightarrow3\rightarrow4\rightarrow5\rightarrow6\rightarrow7\rightarrow8\) \\
\hline
10 & 3, 13 & \(10\rightarrow3\rightarrow6\rightarrow9\rightarrow2\rightarrow5\rightarrow8\rightarrow1\rightarrow4\rightarrow7\rightarrow10\) \\
    & 7 & \(10\rightarrow7\rightarrow4\rightarrow1\rightarrow8\rightarrow5\rightarrow2\rightarrow9\rightarrow6\rightarrow3\rightarrow10\) \\
    & 9 & \(10\rightarrow 9\rightarrow 8\rightarrow 7\rightarrow 6\rightarrow5\rightarrow4\rightarrow3\rightarrow2\rightarrow1\rightarrow10\) \\
    & 11 & \(10\rightarrow1\rightarrow2\rightarrow3\rightarrow4\rightarrow5\rightarrow6\rightarrow7\rightarrow8\rightarrow9\rightarrow10\) \\
\hline
12 & 5 & \(12\rightarrow5\rightarrow10\rightarrow3\rightarrow 8\rightarrow1\rightarrow6\rightarrow11\rightarrow4\rightarrow9\rightarrow2\rightarrow7\rightarrow12\) \\
    & 7 & \(12\rightarrow7\rightarrow2\rightarrow9\rightarrow4\rightarrow11\rightarrow6\rightarrow1\rightarrow8\rightarrow3\rightarrow10\rightarrow5\rightarrow12\) \\
    & 11 & \(12\rightarrow11\rightarrow10\rightarrow9\rightarrow8\rightarrow 7\rightarrow 6\rightarrow5\rightarrow4\rightarrow3\rightarrow2\rightarrow1\rightarrow12\) \\
    & 13 & \(12\rightarrow1\rightarrow2\rightarrow3\rightarrow4\rightarrow5\rightarrow6\rightarrow7\rightarrow8\rightarrow9\rightarrow 10\rightarrow 11\rightarrow12\) \\
\end{longtable}

\bigskip
\subsection{Kolam Patterns as Directed Eulerian Graphs}

The dot connections in a Hridaya Kolam pattern can be modeled as a directed Eulerian graph, where each vertex (dot) has equal in-degree and out-degree. This structure ensures the existence of an Eulerian circuit—a path that visits every directed edge exactly once—enabling the entire pattern to be drawn in a single, uninterrupted motion.

\begin{figure}[H]
    \centering
    \includegraphics[width=0.45\textwidth]{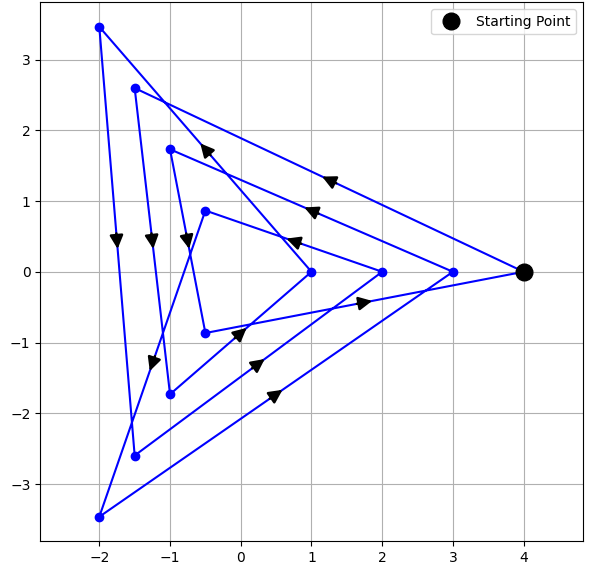}
    \hspace{1em}
    \includegraphics[width=0.45\textwidth]{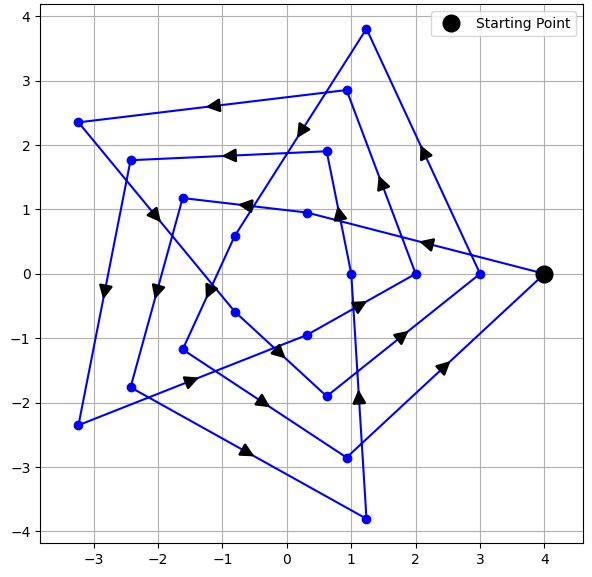}
    \caption{Directed Eulerian graph representations of Hridaya Kolam patterns with even dots: \( m = 4, n = 3 \) (left), \( m = 4, n = 5 \) (right). The Eulerian trail follows all directed edges precisely once, mimicking the continuous drawing of the Kolam.}
\end{figure}

\noindent Each directed edge represents a stroke between adjacent dots, following the Kolam’s orientation. The balance of in-degree and out-degree at every vertex guarantees a continuous path that covers all edges exactly once, without retracing or lifting the pen.

\bigskip

\noindent This graph-theoretic framework enhances the modular arithmetic-based generation of Hridaya Kolam patterns and supports the creation of a wide range of intricate and continuous designs.

\bigskip
\subsection{Algorithm to generate Hridaya Kolam patterns}\label{sec:Algorithm}
To generate the Hridaya Kolam patterns using even-ordered radial dots, we employ the modular sequence that governs the radial distances at each angular step. We fix the number of arms \( n \) (an even number), and generate a base sequence using modular multiplication. This sequence is then repeated around the circle and plotted in polar coordinates to yield closed-loop patterns. The following algorithm provides a step-by-step procedure for computing the polar coordinates required to draw a closed Hridaya Kolam pattern. 
\bigskip
\begin{algorithm}[H]
\caption{Closed-Loop Generation of Even-Ordered-dot Hridaya Kolam}
\begin{algorithmic}[1]
\State \textbf{Input:} Number of radial dots $m$ (even), number of arms $n$ such that $\gcd(m, n) = 1$
\State \textbf{Output:} Set of polar coordinates $(r, \theta)$ to plot a closed Kolam

\State Initialize angle step $\Delta \theta = \frac{2\pi}{n}$
\State Initialize generator sequence $S = [\,]$

\For{$k = 0$ to $m - 1$}
    \If{$k == 0$}
        \State $a_k = m$
    \Else
        \State $a_k = (k \cdot n) \mod m$
    \EndIf
    \State Append $a_k$ to $S$
\EndFor

\State Form extended sequence $S'$ by repeating $S$ exactly $n$ times
\State Initialize coordinate list $C = [\,]$

\For{$i = 0$ to length($S'$)$-1$}
    \State Compute angle $\theta_i = (i \mod n) \cdot \Delta \theta$
    \State Compute radius $r_i = S'[i]$
    \State Append point $(r_i, \theta_i)$ to $C$
\EndFor

\State Append the initial point $C[0]$ at the end of $C$ to close the loop
\State \Return $C$
\end{algorithmic}
\end{algorithm}
 \bigskip

\subsection{Visual Patterns for Generator Sequences}\label{sec:visualpatterns}
The visual patterns presented below depict Hridaya Kolams constructed with an even number of radial dots, where the dot count $m$ is even and the arm count $n$ is chosen such that $\gcd(m, n) = 1$. These patterns illustrate how modular arithmetic sequences can generate symmetric, single-stroke designs inspired by traditional kolam structures. For each fixed even $m$, we vary $n$ over coprime values up to 13 to explore the resulting diversity in geometric arrangements. The designs correspond to the modular sequences previously listed in Table 1. Each subfigure showcases a unique combination of $(m, n)$, revealing rotational symmetries and structural elegance inherent to such configurations.\\\\
All the kolam visuals have been programmatically generated using Python, ensuring consistent styling and algorithmic accuracy. These illustrations serve as a visual verification of the underlying theory used in their construction.

\begin{figure}[H]
    \centering
    \begin{subfigure}[b]{0.3\textwidth}
        \includegraphics[width=\textwidth]{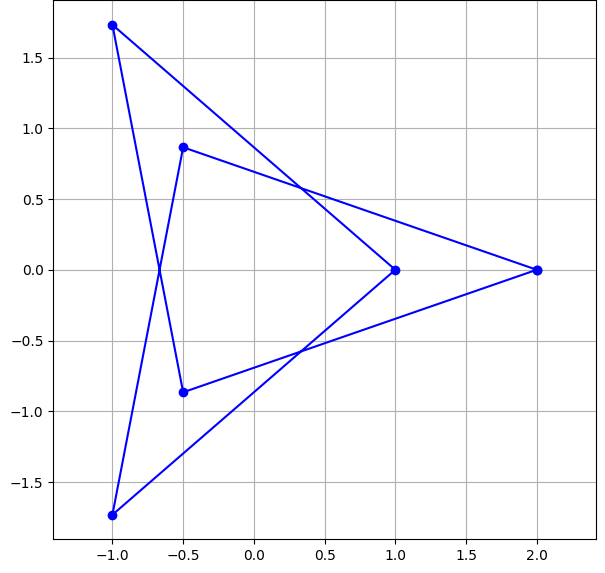}
        \caption{$m=2, n=3$}
    \end{subfigure}\hfill
    \begin{subfigure}[b]{0.3\textwidth}
        \includegraphics[width=\textwidth]{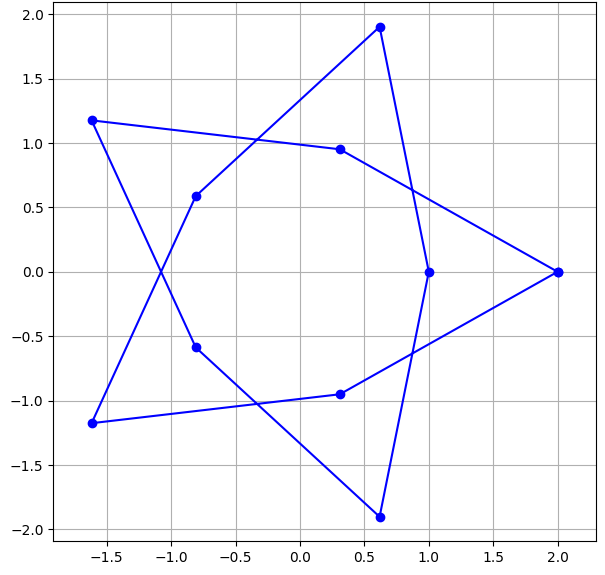}
        \caption{$m=2, n=5$}
    \end{subfigure}\hfill
    \begin{subfigure}[b]{0.3\textwidth}
        \includegraphics[width=\textwidth]{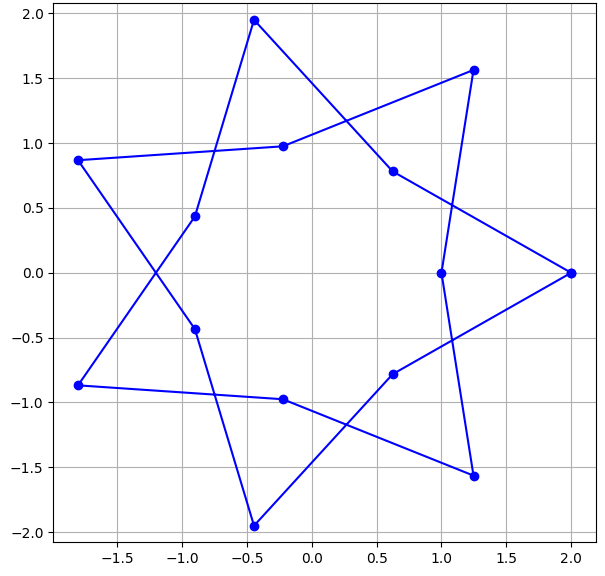}
        \caption{$m=2, n=7$}
    \end{subfigure}
    
    \begin{subfigure}[b]{0.3\textwidth}
        \includegraphics[width=\textwidth]{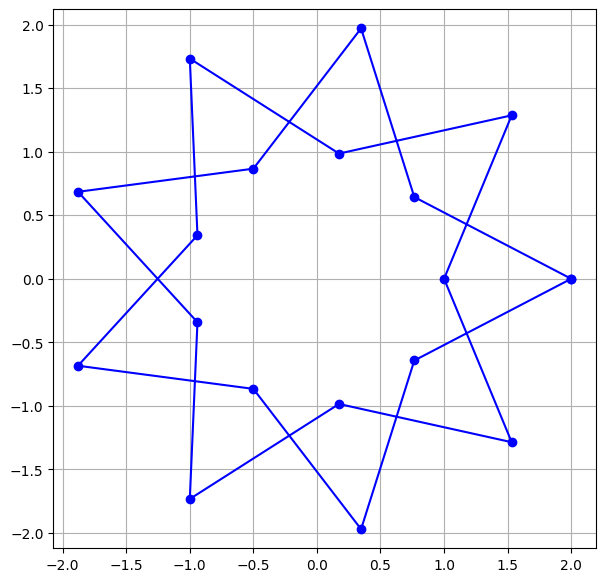}
        \caption{$m=2, n=9$}
    \end{subfigure}\hfill
    \begin{subfigure}[b]{0.3\textwidth}
        \includegraphics[width=\textwidth]{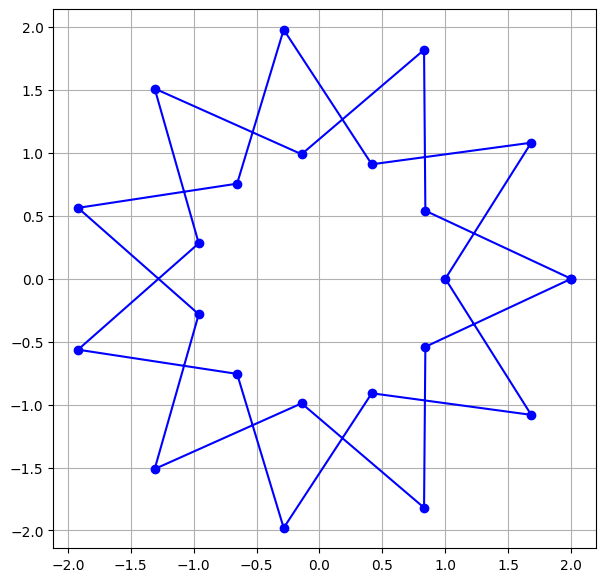}
        \caption{$m=2, n=11$}
    \end{subfigure}\hfill
    \begin{subfigure}[b]{0.3\textwidth}
        \includegraphics[width=\textwidth]{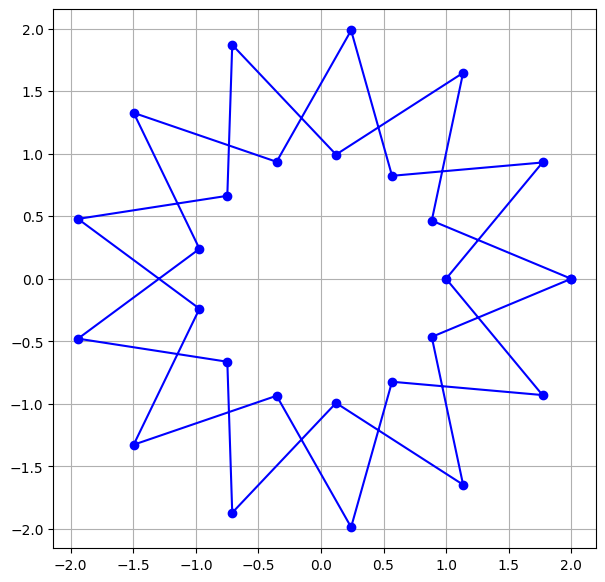}
        \caption{$m=2, n=13$}
    \end{subfigure}
    
    \caption{Patterns for \(m=2\) with various \(n\)}
\end{figure}

\begin{figure}[H]
    \centering
    \begin{subfigure}[b]{0.3\textwidth}
        \includegraphics[width=\textwidth]{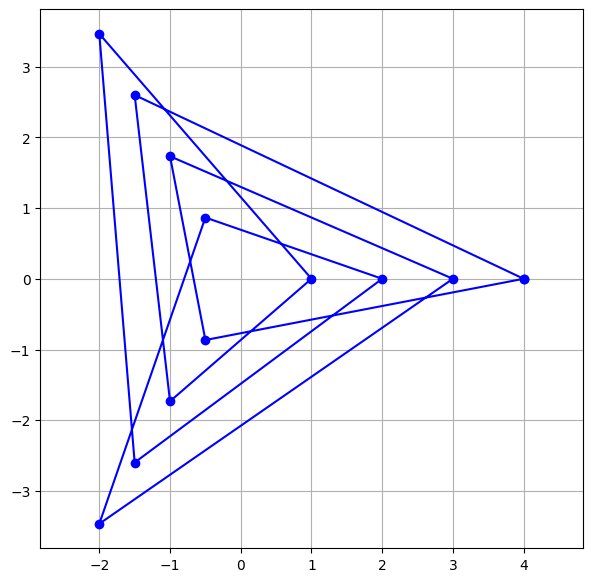}
        \caption{$m=4, n=3$}
    \end{subfigure}\hfill
    \begin{subfigure}[b]{0.3\textwidth}
        \includegraphics[width=\textwidth]{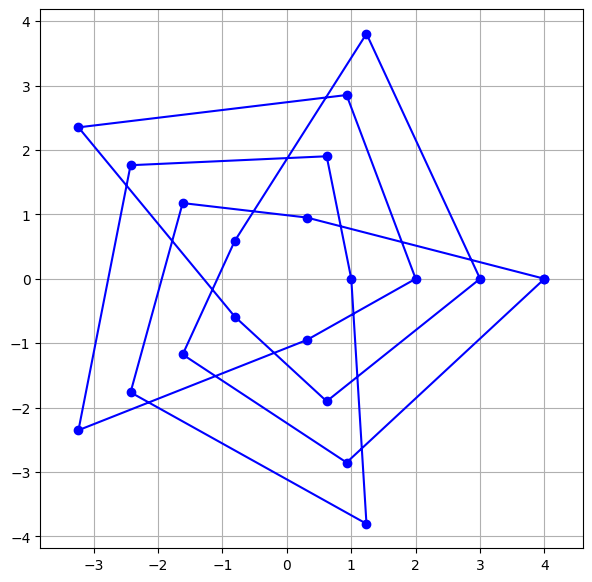}
        \caption{$m=4, n=5$}
    \end{subfigure}\hfill
    \begin{subfigure}[b]{0.3\textwidth}
        \includegraphics[width=\textwidth]{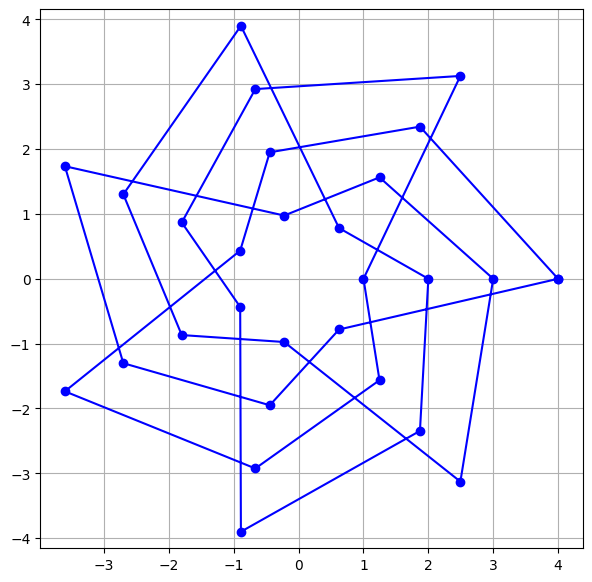}
        \caption{$m=4, n=7$}
    \end{subfigure}
    
    \begin{subfigure}[b]{0.3\textwidth}
        \includegraphics[width=\textwidth]{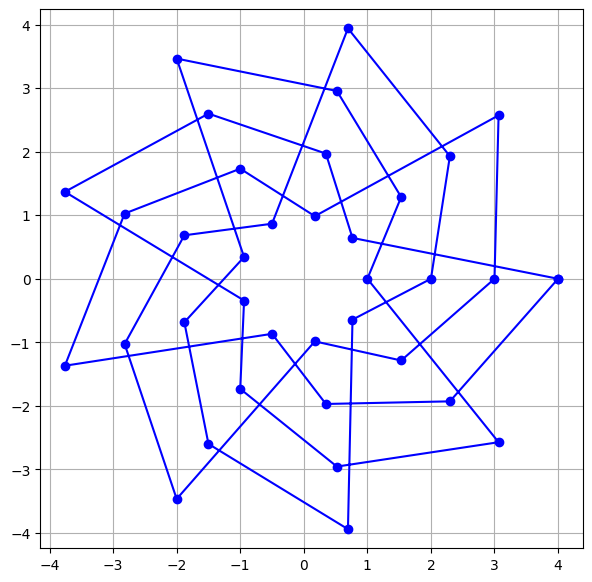}
        \caption{$m=4, n=9$}
    \end{subfigure}\hfill
    \begin{subfigure}[b]{0.3\textwidth}
        \includegraphics[width=\textwidth]{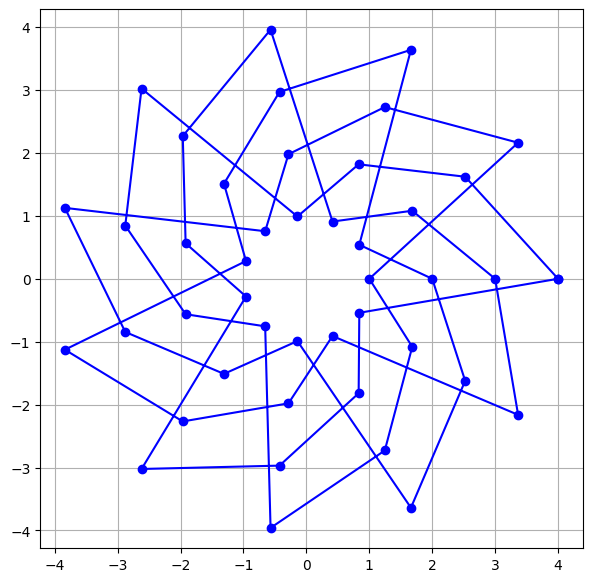}
        \caption{$m=4, n=11$}
    \end{subfigure}\hfill
    \begin{subfigure}[b]{0.3\textwidth}
        \includegraphics[width=\textwidth]{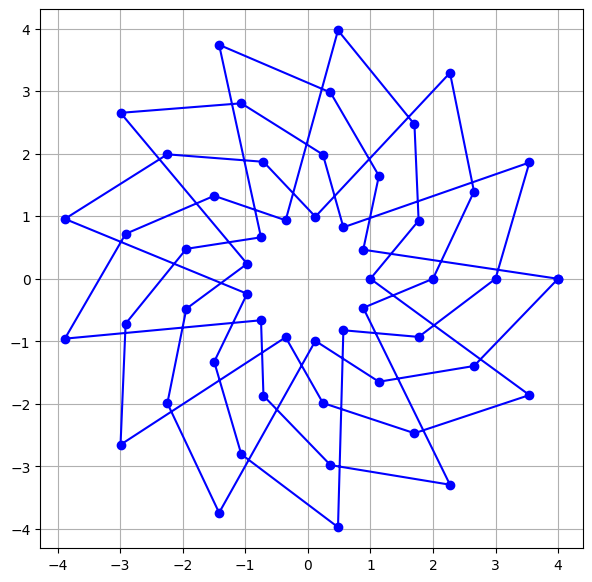}
        \caption{$m=4, n=13$}
    \end{subfigure}
    
    \caption{Patterns for \(m=4\) with various \(n\)}
\end{figure}

\begin{figure}[H]
    \centering
    \begin{subfigure}[b]{0.3\textwidth}
        \includegraphics[width=\textwidth]{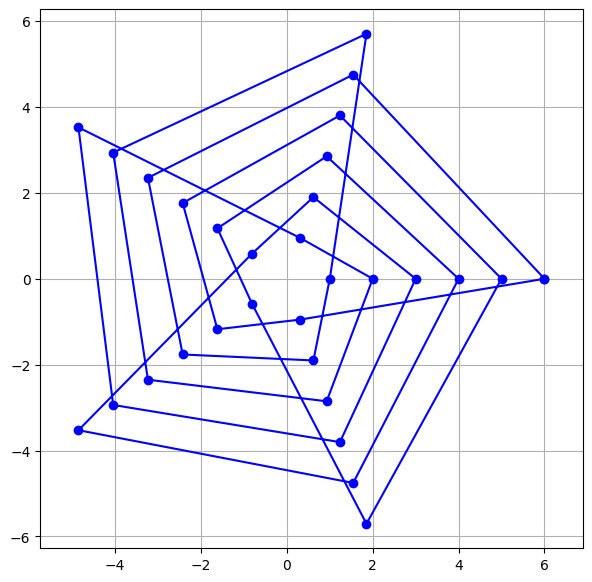}
        \caption{$m=6, n=5$}
    \end{subfigure}\hfill
    \begin{subfigure}[b]{0.3\textwidth}
        \includegraphics[width=\textwidth]{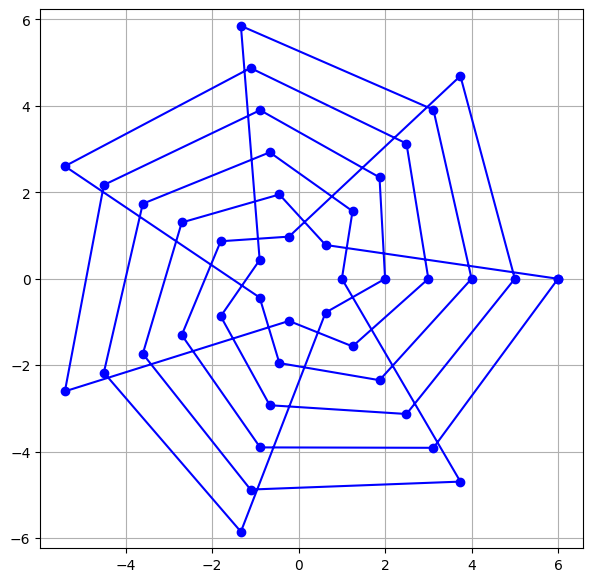}
        \caption{$m=6, n=7$}
    \end{subfigure}\hfill
    \begin{subfigure}[b]{0.3\textwidth}
        \includegraphics[width=\textwidth]{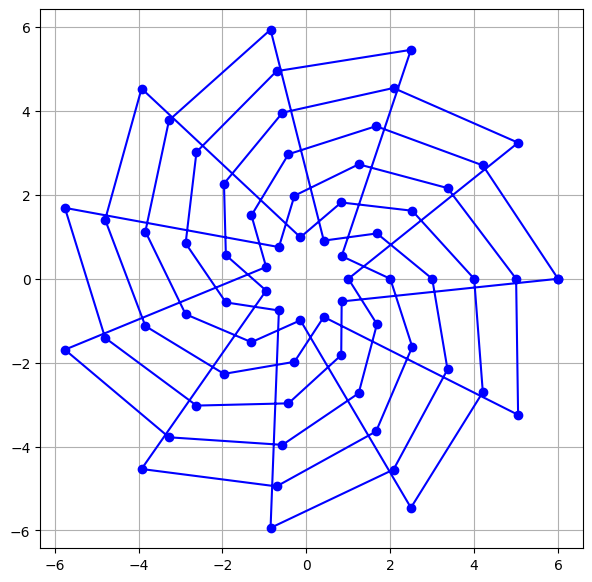}
        \caption{$m=6, n=11$}
    \end{subfigure}
    
    \begin{subfigure}[b]{0.3\textwidth}
        \includegraphics[width=\textwidth]{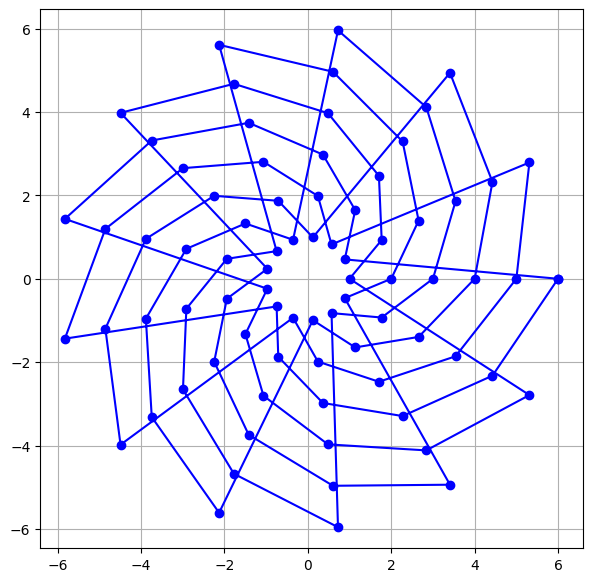}
        \caption{$m=6, n=13$}
    \end{subfigure}
    
    \caption{Patterns for \(m=6\) with various \(n\)}
\end{figure}

\begin{figure}[H]
    \centering
    \begin{subfigure}[b]{0.3\textwidth}
        \includegraphics[width=\textwidth]{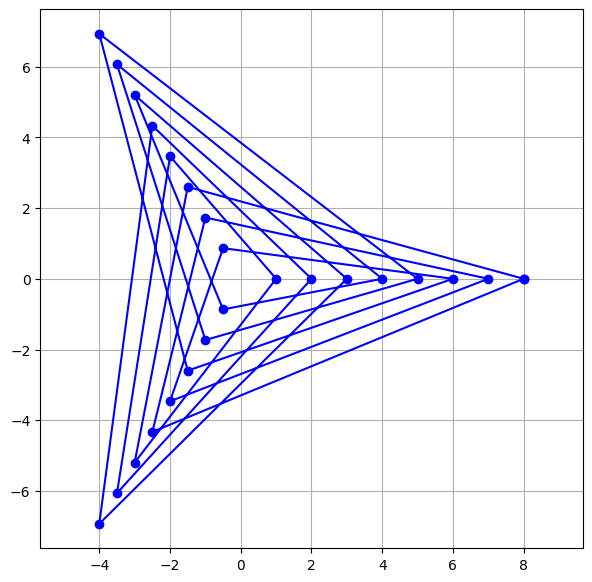}
        \caption{$m=8, n=3$}
    \end{subfigure}\hfill
    \begin{subfigure}[b]{0.3\textwidth}
        \includegraphics[width=\textwidth]{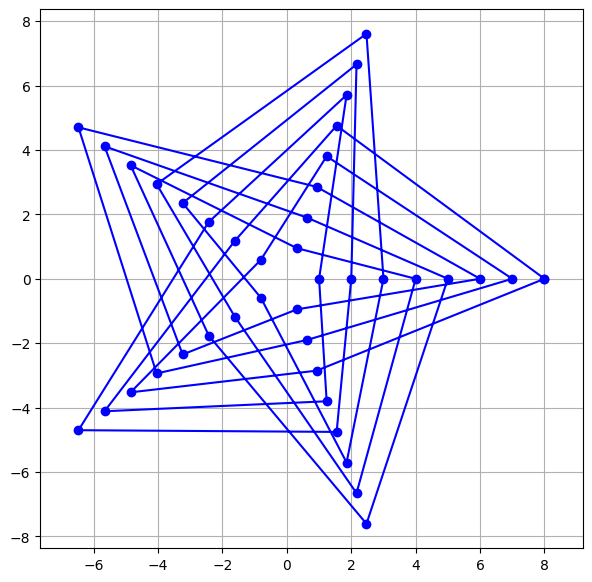}
        \caption{$m=8, n=5$}
    \end{subfigure}\hfill
    \begin{subfigure}[b]{0.3\textwidth}
        \includegraphics[width=\textwidth]{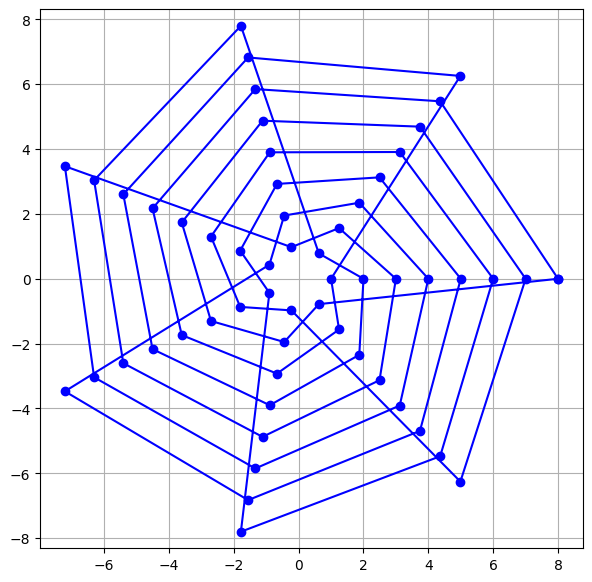}
        \caption{$m=8, n=7$}
    \end{subfigure}
    
    \begin{subfigure}[b]{0.3\textwidth}
        \includegraphics[width=\textwidth]{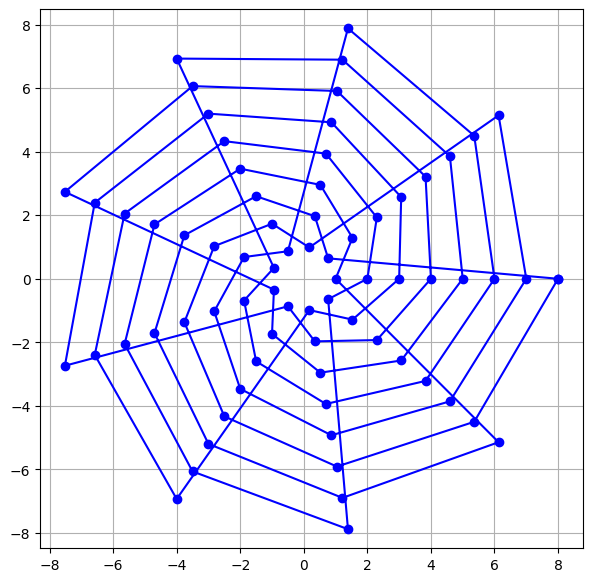}
        \caption{$m=8, n=9$}
    \end{subfigure}\hfill
    \begin{subfigure}[b]{0.3\textwidth}
        \includegraphics[width=\textwidth]{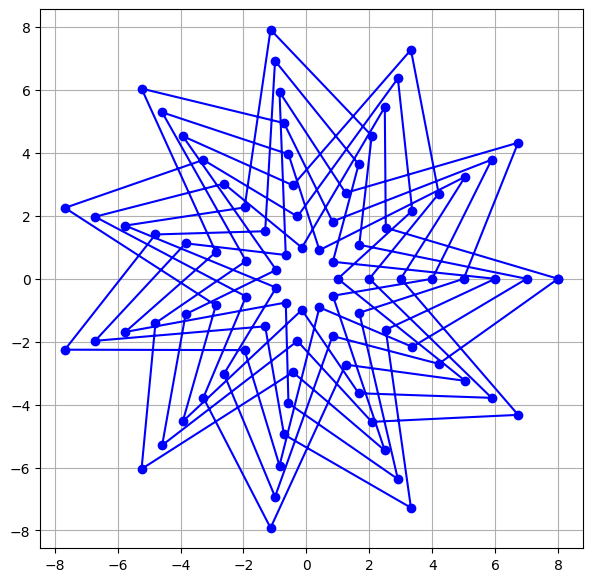}
        \caption{$m=8, n=11$}
    \end{subfigure}\hfill
    \begin{subfigure}[b]{0.3\textwidth}
        \includegraphics[width=\textwidth]{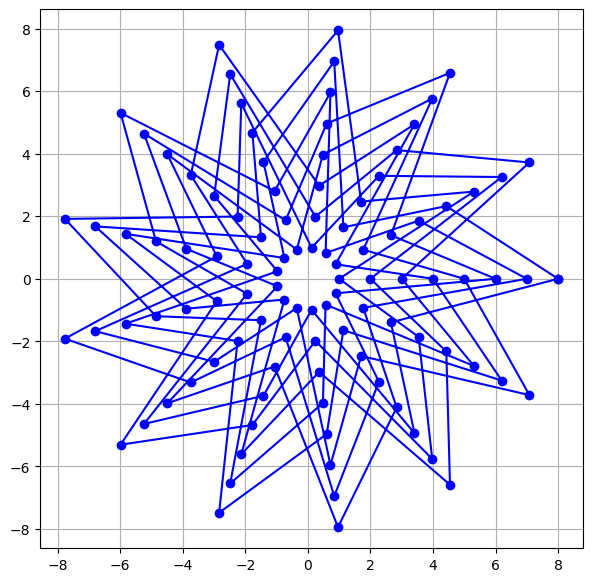}
        \caption{$m=8, n=13$}
    \end{subfigure}
    
    \caption{Patterns for \(m=8\) with various \(n\)}
\end{figure}

\begin{figure}[H]
    \centering
    \begin{subfigure}[b]{0.3\textwidth}
        \includegraphics[width=\textwidth]{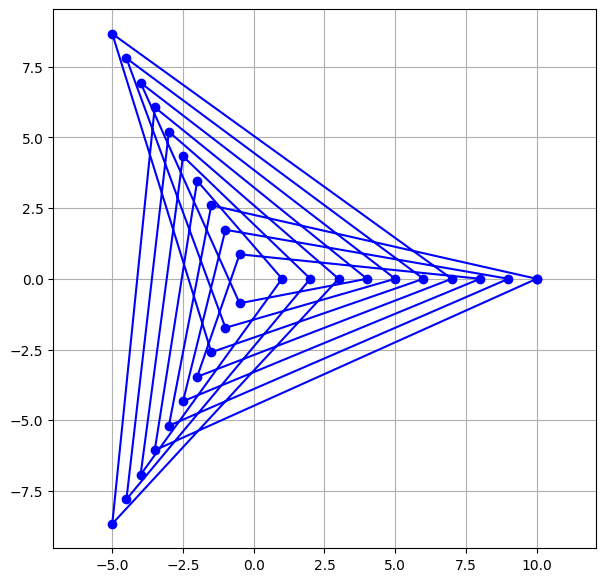}
        \caption{$m=10, n=3$}
    \end{subfigure}\hfill
    \begin{subfigure}[b]{0.3\textwidth}
        \includegraphics[width=\textwidth]{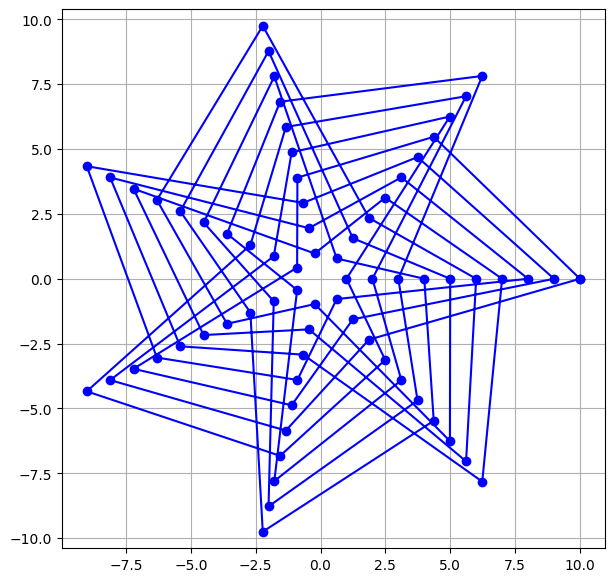}
        \caption{$m=10, n=7$}
    \end{subfigure}\hfill
    \begin{subfigure}[b]{0.3\textwidth}
        \includegraphics[width=\textwidth]{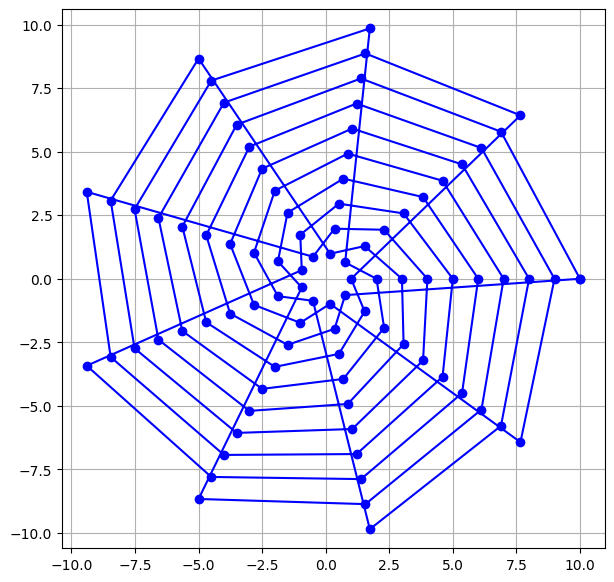}
        \caption{$m=10, n=9$}
    \end{subfigure}

    \centering
    \begin{subfigure}[b]{0.3\textwidth}
        \includegraphics[width=\textwidth]{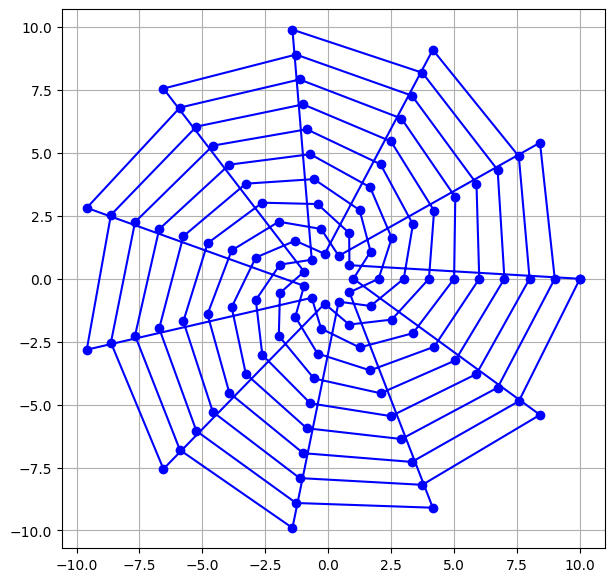}
        \caption{$m=10, n=11$}
    \end{subfigure}
    \hspace{0.1\textwidth} 
    \begin{subfigure}[b]{0.3\textwidth}
        \includegraphics[width=\textwidth]{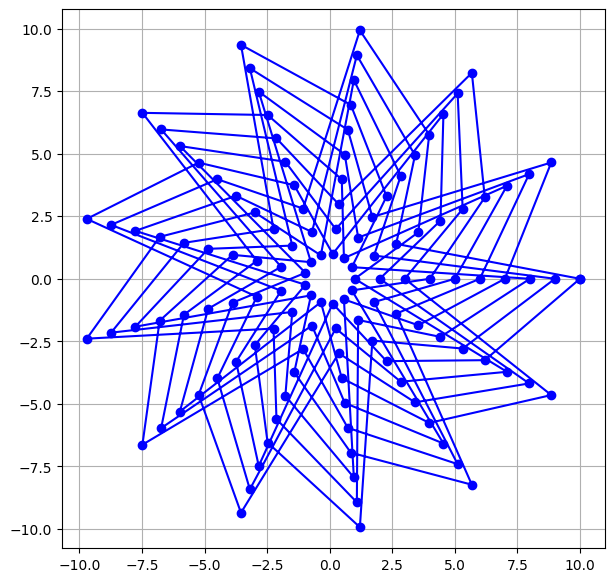}
        \caption{$m=10, n=13$}
    \end{subfigure}

    \caption{Patterns for \(m=10\) with various \(n\)}
\end{figure}

\begin{figure}[H]
    \centering
    \begin{subfigure}[b]{0.3\textwidth}
        \includegraphics[width=\textwidth]{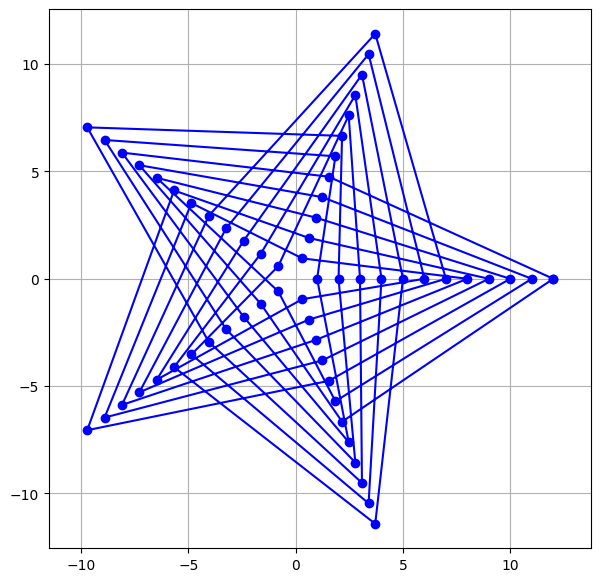}
        \caption{$m=12, n=5$}
    \end{subfigure}\hfill
    \begin{subfigure}[b]{0.3\textwidth}
        \includegraphics[width=\textwidth]{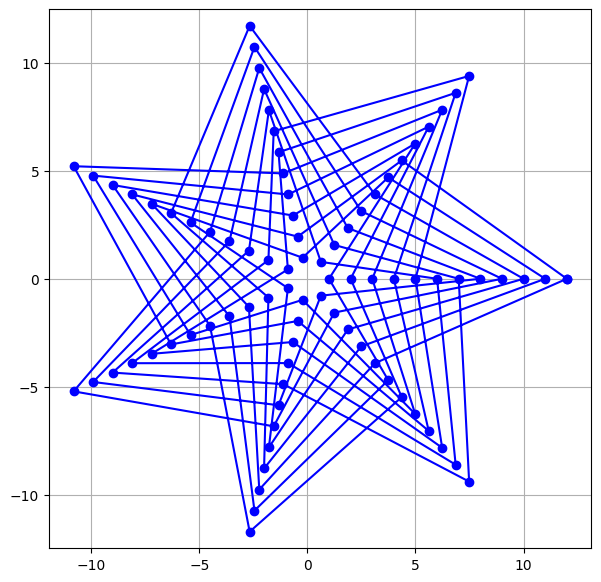}
        \caption{$m=12, n=7$}
    \end{subfigure}\hfill
    \begin{subfigure}[b]{0.3\textwidth}
        \includegraphics[width=\textwidth]{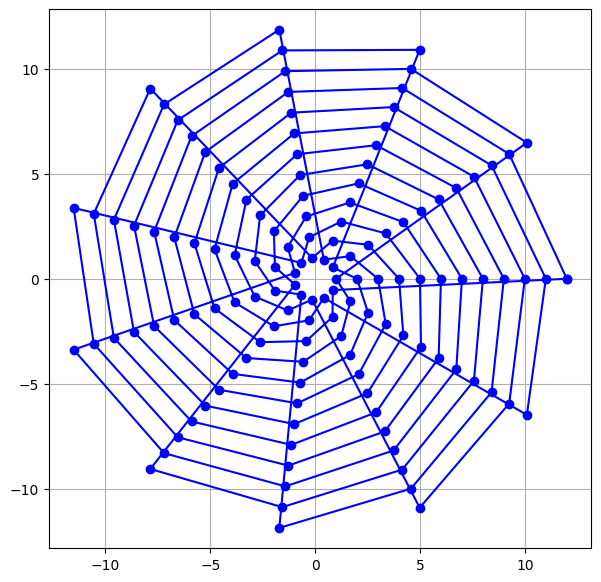}
        \caption{$m=12, n=11$}
    \end{subfigure}

    \begin{subfigure}[b]{0.3\textwidth}
        \includegraphics[width=\textwidth]{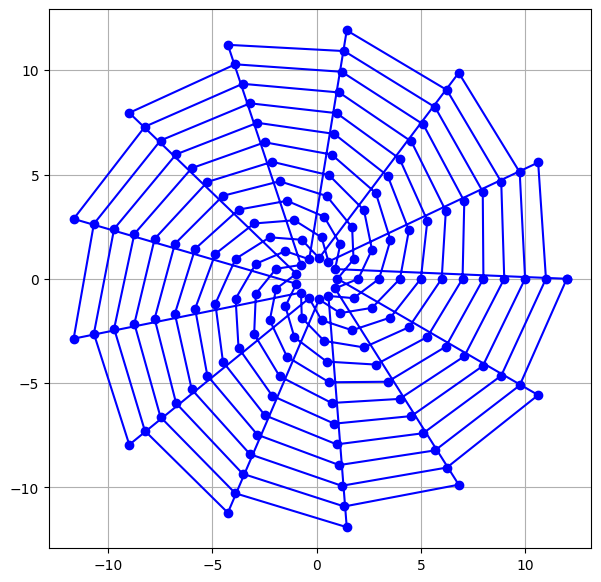}
        \caption{$m=12, n=13$}
    \end{subfigure}
    
    \caption{Patterns for \(m=12\) with various \(n\)}
\end{figure}


\subsection{Extension to Larger Arm Counts with Even Dots}

To explore more complex Hridaya Kolam structures, we fix the number of radial dots at \( m = 20 \) and vary the arm count \( n \) over values coprime to \( m \), ensuring \( \gcd(m,n) = 1 \). Fixing \( m = 20 \) strikes a balance—keeping the dot structure visually clear while still allowing intricate designs. Larger values of both \( m \) and \( n \) tend to saturate the image, making the dot layout hard to perceive. Below are visualizations for selected values of \( n \): \( 7, 13, 19, 23, 27, 91 \). The final case with \( n = 91 \) approaches a visual limit, nearly filling the circular space and producing a continuous-looking pattern.

\begin{figure}[H]
    \centering
    \begin{subfigure}[b]{0.3\textwidth}
        \includegraphics[width=\textwidth]{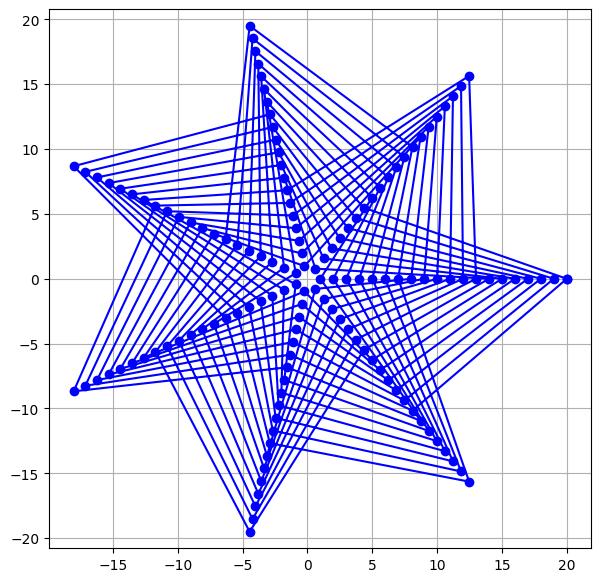}
        \caption{Pattern for \(n=7\)}
    \end{subfigure}\hfill
    \begin{subfigure}[b]{0.3\textwidth}
        \includegraphics[width=\textwidth]{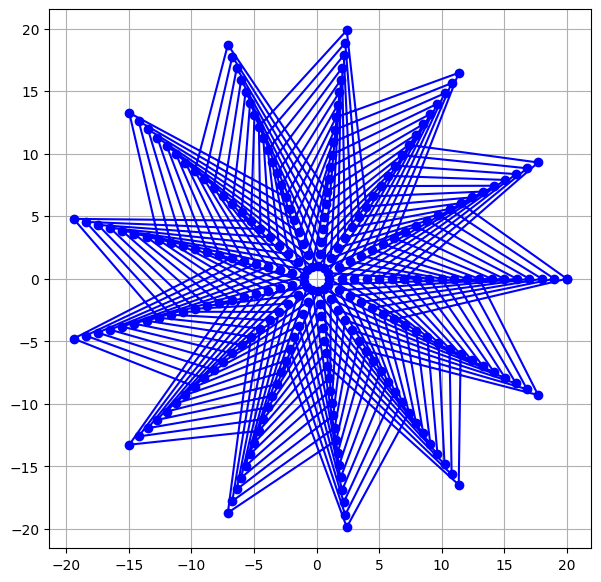}
        \caption{Pattern for \(n=13\)}
    \end{subfigure}\hfill
    \begin{subfigure}[b]{0.3\textwidth}
        \includegraphics[width=\textwidth]{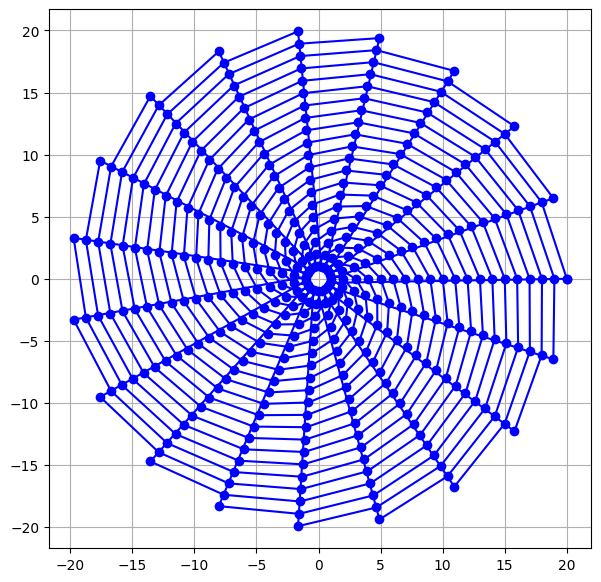}
        \caption{Pattern for \(n=19\)}
    \end{subfigure}
    
    \vspace{0.5cm}
    
    \begin{subfigure}[b]{0.3\textwidth}
        \includegraphics[width=\textwidth]{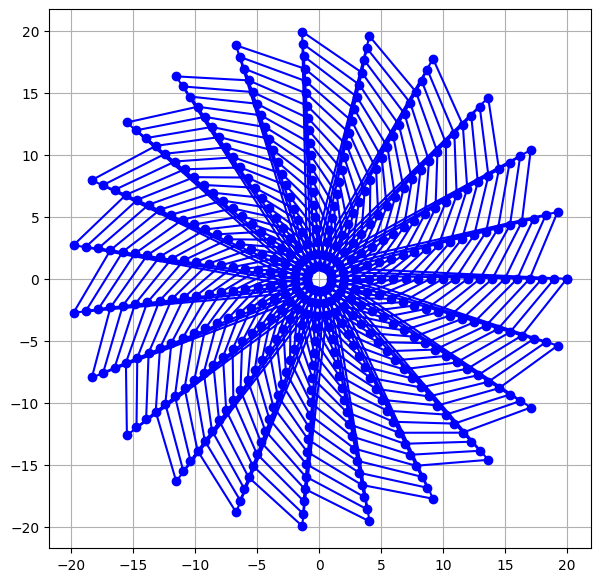}
        \caption{Pattern for \(n=23\)}
    \end{subfigure}\hfill
    \begin{subfigure}[b]{0.3\textwidth}
        \includegraphics[width=\textwidth]{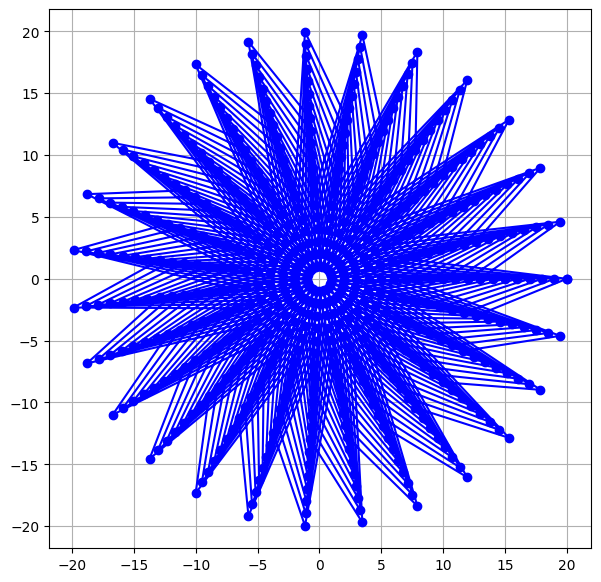}
        \caption{Pattern for \(n=27\)}
    \end{subfigure}\hfill
    \begin{subfigure}[b]{0.3\textwidth}
        \includegraphics[width=\textwidth]{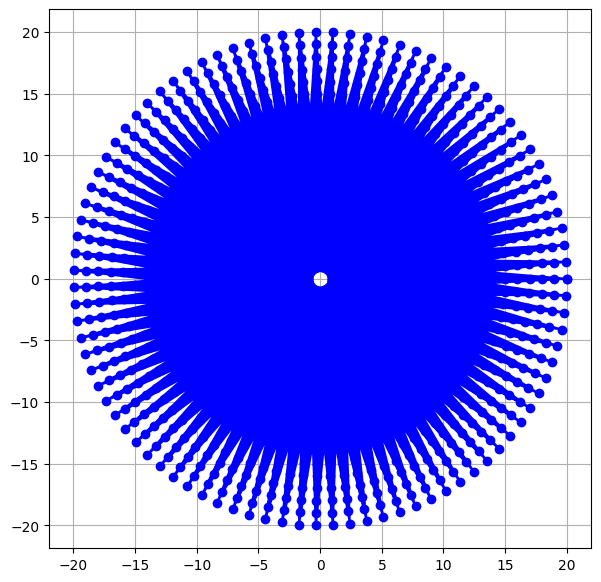}
        \caption{Pattern for \(n=91\)}
    \end{subfigure}
    
    \caption{Hridaya Kolam visualizations with fixed \(m = 20\) and increasing \(n\), each coprime to 20.}
    \label{fig:extended_even_hridaya}
\end{figure}

\subsection{Visualization of Kolam Patterns Using Circular Arcs}

Continuing with the dot-arm structure for Hridaya Kolams with even \( m \), we now explore different styles of connecting the dots. Traditionally, \emph{bulging outward} (convex) arcs are used, curving away from the center to produce smooth, elegant loops. To broaden the visual vocabulary, we also consider \emph{bulging inward} (concave) arcs, which curve toward the center and offer a distinct aesthetic, though less common in traditional designs. Additionally, we include the simpler case of \emph{straight-line} connections, previously used, as a baseline for comparison. The figure below shows all three styles for \((m, n) = (4, 5)\).

\begin{figure}[H]
    \centering
    \begin{subfigure}[b]{0.3\textwidth}
        \centering
        \includegraphics[width=\linewidth]{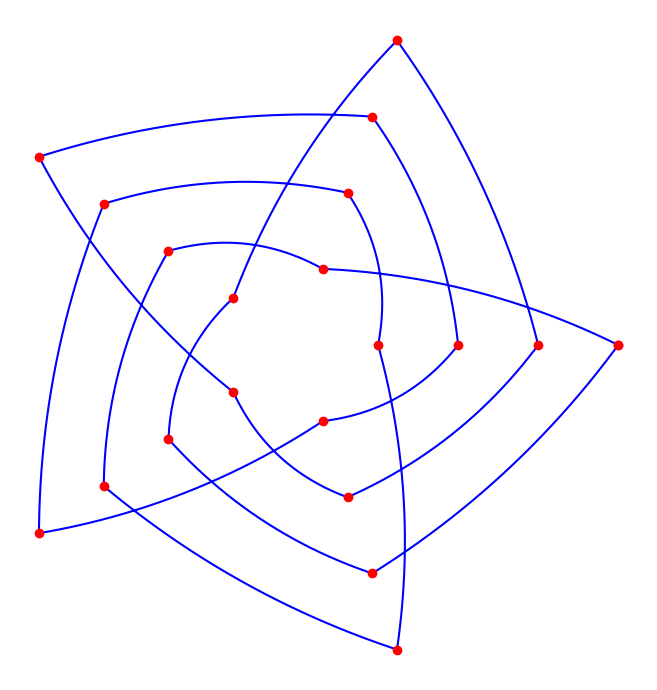}
        \caption{Bulging outward arcs}
    \end{subfigure}
    \hfill
    \begin{subfigure}[b]{0.3\textwidth}
        \centering
        \includegraphics[width=\linewidth]{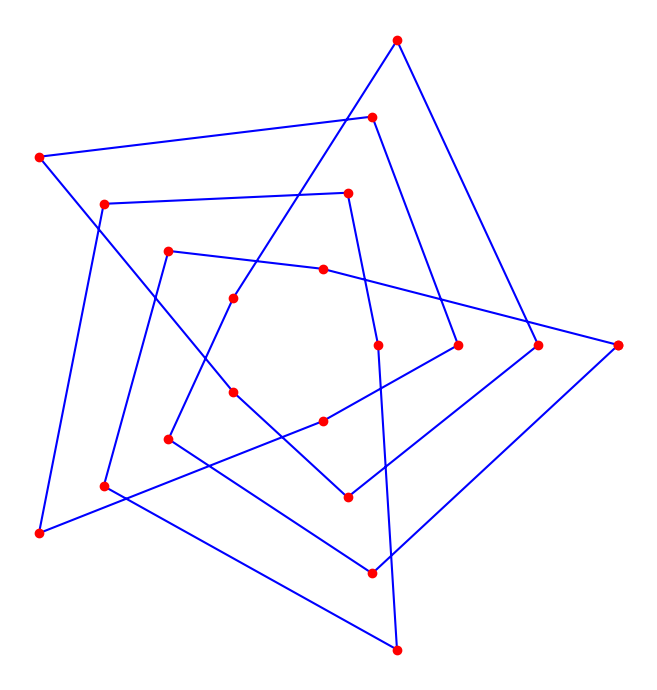}
        \caption{Straight lines}
    \end{subfigure}
    \hfill
    \begin{subfigure}[b]{0.3\textwidth}
        \centering
        \includegraphics[width=\linewidth]{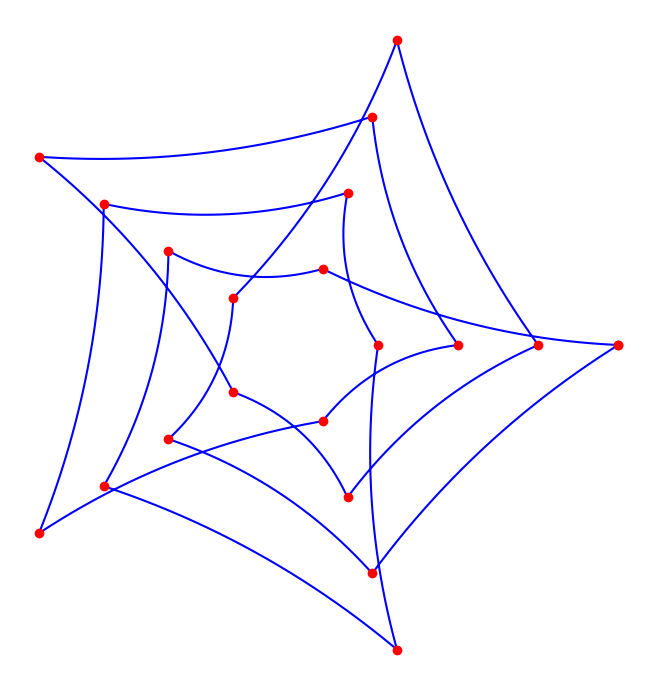}
        \caption{Bulging inward arcs}
    \end{subfigure}
    
    \caption{Hridaya Kolam for \((m,n) = (4,5)\), using three connection styles: outward arcs, straight lines, and inward arcs.}
    \label{fig:kolam_arcs}
\end{figure}

\section{Application in the Carpet Sector}\label{sec:application}

We applied the Even-Dotted Hridaya Kolam framework to create intricate designs using Python for pattern generation and CorelDRAW for visualization. Traditionally, Kolams are ephemeral artworks made on the ground with rice flour or chalk. Here, we reinterpret them as permanent, functional textile designs. By selecting appropriate color schemes for the subregions of the patterns, we transformed these mathematical motifs into woven products. Several Daries were produced at the Indian Institute of Carpet Technology, demonstrating a thoughtful blend of traditional artistry and mathematical aesthetics.

\begin{figure}[H]
\centering
\begin{subfigure}[b]{0.27\textwidth}
    \centering
    \includegraphics[width=\linewidth]{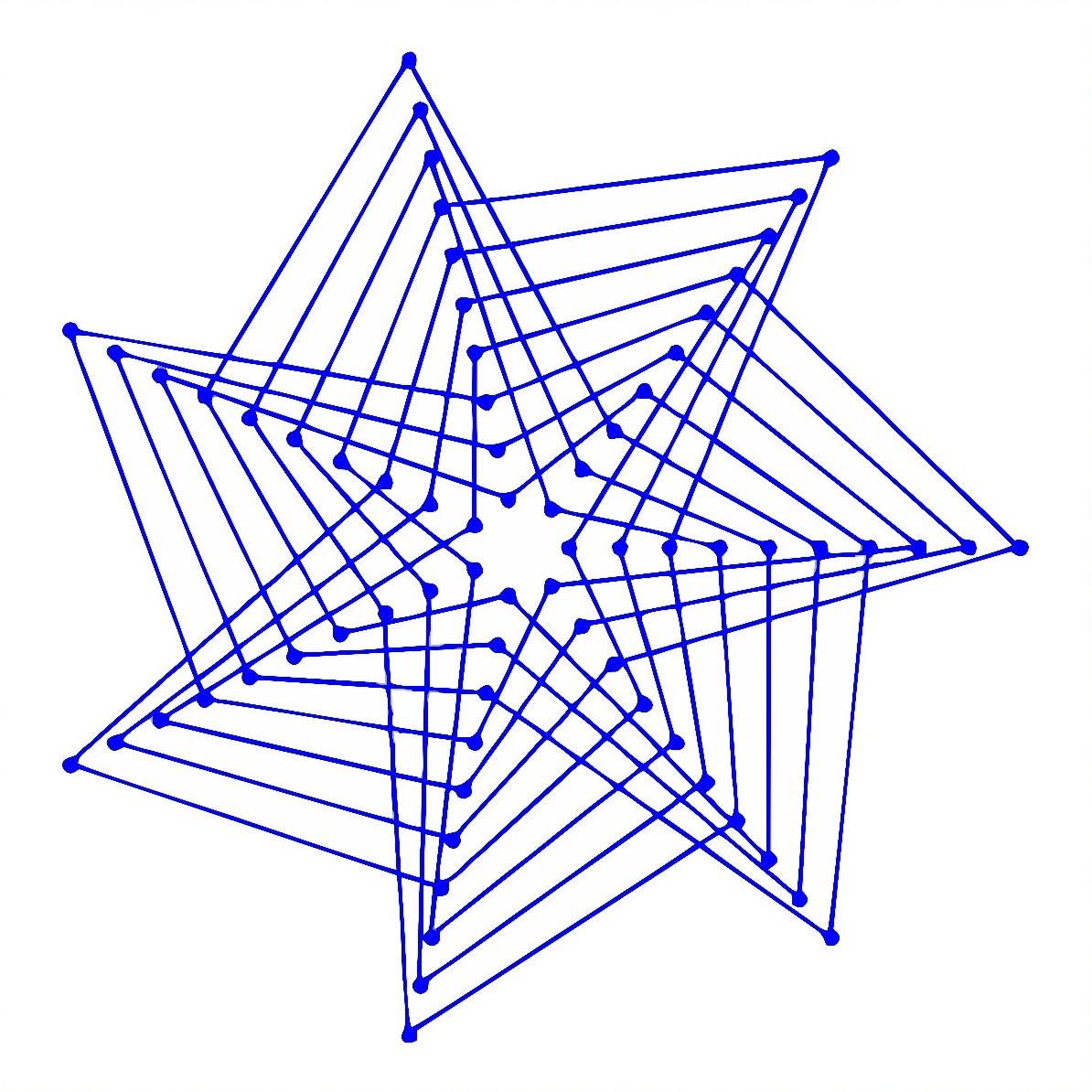}
    \caption{Design (CorelDRAW)}
    \label{fig:coreldraw_plot}
\end{subfigure}
\hfill
\begin{subfigure}[b]{0.3\textwidth}
    \centering
    \includegraphics[width=\linewidth]{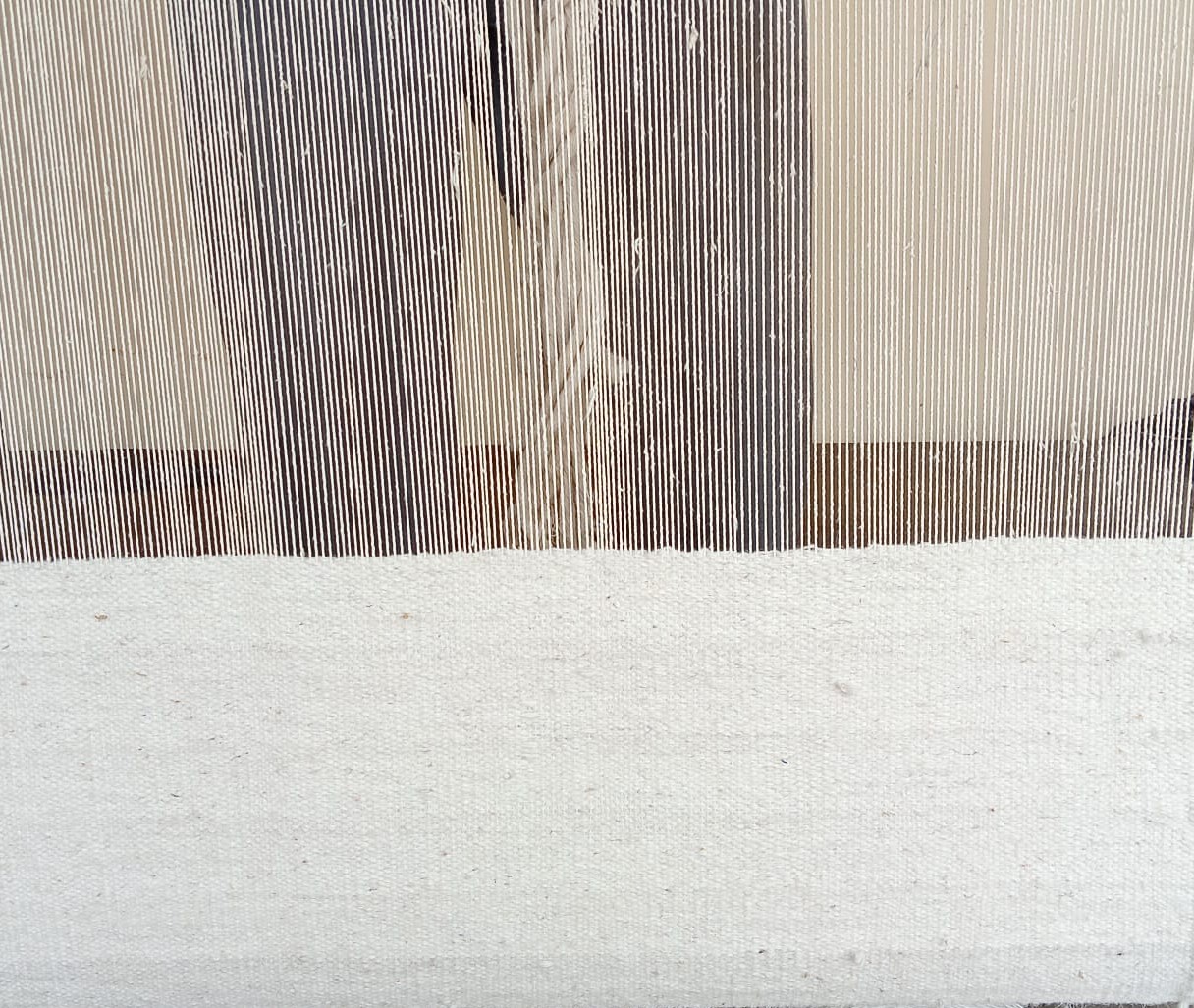}
    \caption{During Dari weaving}
    \label{fig:during_manufacturing}
\end{subfigure}
\hfill
\begin{subfigure}[b]{0.31\textwidth}
    \centering
    \includegraphics[width=\linewidth]{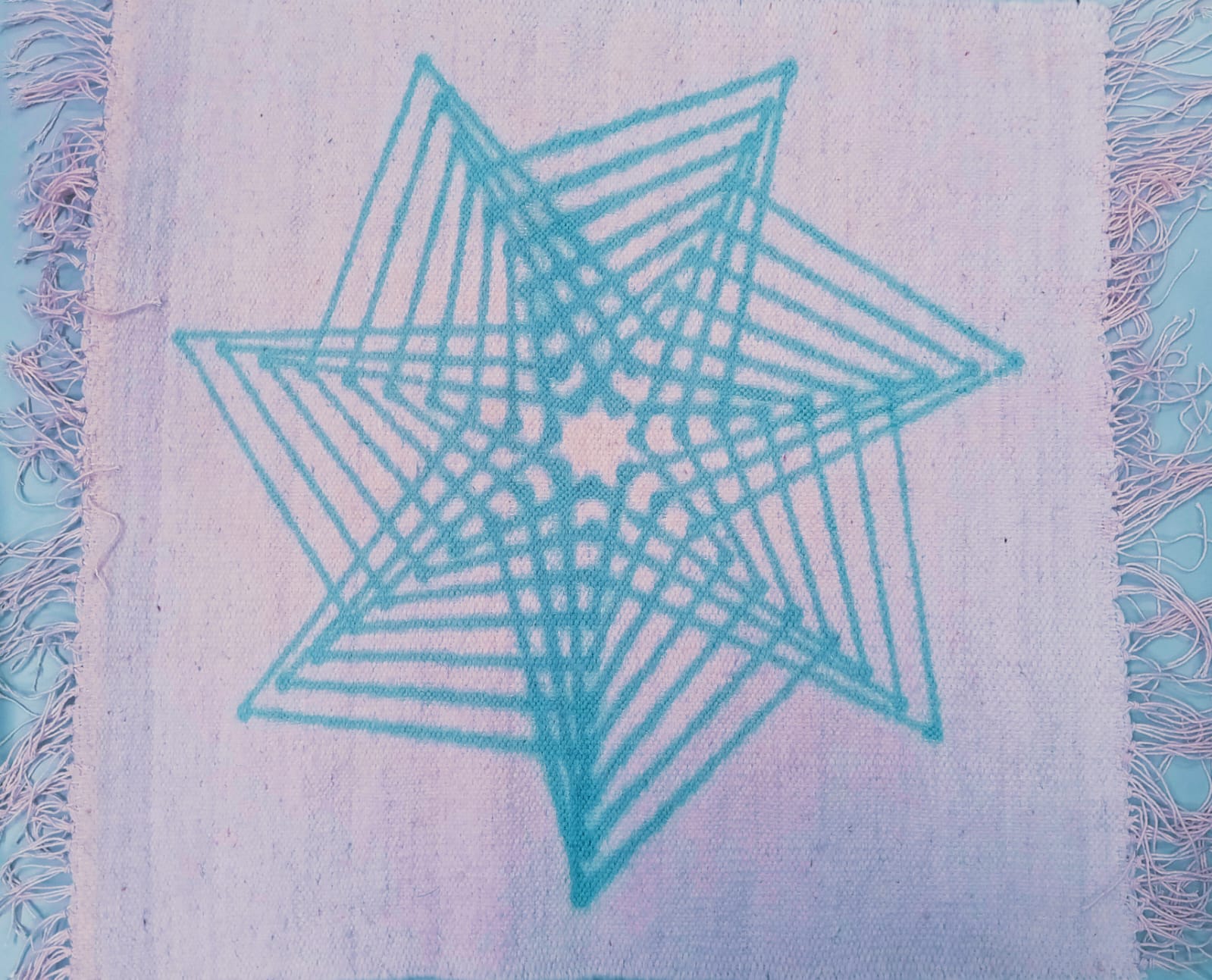}
    \caption{Finished woven Dari}
    \label{fig:final_dari}
\end{subfigure}

\caption{Kolam pattern transformation for $(m,n) = (10,7)$. The initial design created in Python (See Figure~10(b)) which was refined in CorelDRAW (left), prepared for weaving (middle), and finally rendered as a finished Dari (right).}
\label{fig:kolam_to_dari}
\end{figure}

\section{Conclusions} \label{sec:conclusion}

This article explores the connection between Hridaya Kolam designs and number theory, extending pattern generation from odd to even dot configurations using coprime-based modular sequences. By uniting mathematical structure with traditional aesthetics, we generate intricate, symmetric patterns with potential for practical and artistic applications. The approach highlights how culturally rooted designs can be formalized and expanded through algorithmic methods. Future directions include automating Kolam generation, applying these techniques to other motifs, and integrating machine learning for pattern synthesis and recognition—bridging mathematics, art, and technology in meaningful ways.\\

\noindent \textbf{Data availability statement}: Not applicable.\\

\end{document}